\documentclass[12pt,fullpage]{article}

\usepackage{amsfonts,graphicx,amsmath,xr,epsfig,comment}
\input epsf

\def\R{\mathbb{R}}
\def\N{\mathbb{N}}
\def\epsilon{\varepsilon}
\def\trait (#1) (#2) (#3){\vrule width #1pt height #2pt depth #3pt}
\def\fin{\hfill\trait (0.1) (5) (0) \trait (5) (0.1) (0) \kern-5pt \trait (5) (5) (-4.9) \trait (0.1) (5) (0)}

\newcommand{\SE}{\setcounter{equation}{0} \section}
\newcommand{\be}{\begin{equation}}
\newcommand{\ee}{\end{equation}}
\newcommand{\baa}{\begin{array}}
\newcommand{\eaa}{\end{array}}
\newcommand{\ba}{\begin{eqnarray}}
\newcommand{\ea}{\end{eqnarray}}

\newtheorem{theo}{\bf Theorem}[section]

\newtheorem{cor}[theo]{\bf Corollary}

\newtheorem{rem}[theo]{\bf Remark}

\begin{document}
\date{}
\title{\bf{Fast propagation for KPP equations with slowly decaying initial conditions}}
\author{Fran\c cois Hamel$^{\hbox{\small{ a}}}$$\,$ and Lionel Roques$^{\hbox{\small{ b}}}$\thanks{The authors are supported by the French ``Agence Nationale de la Recherche" within the projects ColonSGS, PREFERED, and URTICLIM (second author). The first author is also indebted to the Alexander von~Humboldt Foundation for its support.}\\
\\
\footnotesize{$^{\hbox{a }}$Aix-Marseille Universit\'e, LATP, Facult\'e des Sciences et Techniques}\\
\footnotesize{Avenue Escadrille Normandie-Niemen, F-13397 Marseille Cedex 20, France}\\
\footnotesize{\& Helmholtz Zentrum M\"unchen, Institut f\"ur Biomathematik und Biometrie}\\
\footnotesize{Ingolst\"adter Landstrasse 1, D-85764 Neuherberg, Germany}\\
\footnotesize{E-mail address: francois.hamel@univ-cezanne.fr (corresponding author)}\\
\footnotesize{$^{\hbox{b }}$UR 546 Biostatistique et Processus Spatiaux, INRA, F-84000 Avignon, France}\\
\footnotesize{E-mail address: lionel.roques@avignon.inra.fr}}
\maketitle

\begin{abstract}
This paper is devoted to the analysis of the large-time behavior of solutions of one-dimensional Fisher-KPP reaction-diffusion equations. The initial conditions are assumed to be globally front-like and to decay at infinity towards the unstable steady state more slowly than any exponentially decaying function. We prove that all level sets of the solutions move infinitely fast as time goes to infinity. The locations of the level sets are expressed in terms of the decay of the initial condition. Furthermore, the spatial profiles of the solutions become asymptotically uniformly flat at large time. This paper contains the first systematic study of the large-time behavior of solutions of KPP equations with slowly decaying initial conditions. Our results are in sharp contrast with the well-studied case of exponentially bounded initial conditions.
\end{abstract}



\SE{Introduction}\label{intro}

In this paper, we study the large-time behavior of the solutions of the Cauchy problem for reaction-diffusion equations
\be\label{cauchy}\left\{\baa{l}
u_t=u_{xx}+f(u),\ \ t>0,\ \ x\in\R,\vspace{5pt}\\
u(0,x)=u_0(x),\ \ x\in\R\eaa\right.
\ee
with Fisher-KPP nonlinearities $f$ and with slowly decaying initial conditions~$u_0$. The functions~$f$ are assumed to have two zeroes: the one is unstable and the other one is stable. We prove that, provided that $u_0(x)$ decays sufficiently slowly to the unstable steady state as $x\to+\infty$, then the central part of the solution $u$ moves to the right with infinite speed at large time, in a sense to be described below. Furthermore, the initial condition $u_0$ can be chosen so that the location of the solution~$u$ be asymptotically larger than any prescribed real-valued function.\par
Many papers have been concerned with the large-time behavior of the solutions of equation~(\ref{cauchy}) or more general reaction-diffusion equations with exponentially decaying initial conditions, leading to finite propagation speeds. The results of this paper shed a new light on this problem and they are all the more important as they are already valid for the simplest reaction-diffusion model, that is problem (\ref{cauchy}) in a one-dimensional homogeneous medium. Similar results would hold for more general equations, but we choose to restrict the analysis to the one-dimensional case for the sake of simplicity of the presentation.\par
Let us now make the assumptions more precise. The unknown quantity $u$ typically stands for a density and ranges in the interval $[0,1]$. The given function $f:[0,1]\to\R$ is of class $C^1$ and it is assumed to be of the Fisher-KPP type~\cite{f,kpp}, that is
\be\label{fkpp}
f(0)=f(1)=0,\quad 0<f(s)\le f'(0)\,s\ \hbox{ for all }s\in(0,1).
\ee
The last property means that the growth rate $f(s)/s$ is maximal at $s=0$. Furthermore, throughout the paper, we assume that there exist $\delta>0$, $s_0\in(0,1)$ and $M\ge 0$ such that
$$f(s)\ge f'(0)s-M\,s^{1+\delta}\ \hbox{ for all }s\in[0,s_0].$$
A particular class of such functions $f$ is that of $C^{1,\delta}([0,1])$ concave functions $f$, which are positive in~$(0,1)$ and vanish at $0$ and $1$.\par
The initial condition $u_0:\R\to[0,1]$ is assumed to be uniformly continuous and asymptotically front-like, that~is
\be\label{frontlike}
u_0>0\hbox{ in }\R,\ \ \liminf_{x\to-\infty}u_0(x)>0\ \hbox{ and }\ \lim_{x\to+\infty}u_0(x)=0.
\ee
However, it is worth noticing that $u_0$ is not assumed to be nonincreasing. The maximum principle implies immediately that
$$0<u(t,x)<1\ \hbox{ for all }t>0\ \hbox{ and }x\in\R.$$
Furthermore, the function $u_0$ is assumed to decay more slowly than any exponentially decaying function as $x\to+\infty$, in the sense that
\be\label{slow}
\forall\,\epsilon>0,\ \exists\,x_{\epsilon}\in\R,\quad u_0(x)\ge e^{-\epsilon x}\hbox{ in }[x_{\epsilon},+\infty),
\ee
or, equivalently, $u_0(x)\,e^{\epsilon x}\to+\infty$ as $x\to+\infty$ for all $\epsilon>0$. Condition~(\ref{slow}) is fulfilled in particular if $u_0$ is of class $C^1(\xi_0,+\infty)$ for some $\xi_0\in\R$ and if $u'_0(x)=o(u_0(x))$ as $x\to+\infty$.\par
Before going further on, let us list at least four typical classes of initial conditions~$u_0$ satisfying~(\ref{slow}): 1)~the functions which are logarithmically sublinear as $x\to+\infty$, that is
\be\label{ex1}
u_0(x)\sim C\,e^{-\alpha\,x/\ln x}\hbox{ as }x\to+\infty
\ee
with $\alpha,\ C>0$; 2)~the functions which are logarithmically power-like and sublinear as $x\to+\infty$, that is
\be\label{ex2}
u_0(x)\sim C\,e^{-\beta\,x^{\alpha}}\hbox{ as }x\to+\infty
\ee
with $\alpha\in(0,1)$ and $\beta,\ C>0$; 3)~the functions which decay algebraically as $x\to+\infty$, that~is
\be\label{ex3}
u_0(x)\sim C\,x^{-\alpha}\hbox{ as }x\to+\infty
\ee
with $\alpha,\ C>0$; 4)~the functions which decay as a negative power of $\ln(x)$ as $x\to+\infty$, that is
\be\label{ex4}
u_0(x)\sim C\,(\ln x)^{-\alpha}\hbox{ as }x\to+\infty
\ee
with $\alpha,\ C>0$.\par
This type of assumption (\ref{slow}) is in contrast with the large literature on pro\-blem~(\ref{cauchy}), where the initial conditions $u_0$ are assumed to be exponentially bounded as $x\to+\infty$ and where the solutions $u$ converge in some sense to some traveling fronts with finite speed as $t\to+\infty$. It is indeed well-known that the equation $u_t=u_{xx}+f(u)$ admits a family of traveling fronts $u(t,x)=\varphi_c(x-ct)$ for all speeds $c\ge c^*=2\sqrt{f'(0)}$, where, for each speed $c\in[c^*,+\infty)$, the function $\varphi_c:\R\to(0,1)$ satisfies
\be\label{varphic}
\varphi_c''+c\varphi_c'+f(\varphi_c)=0\hbox{ in }\R,\ \ \varphi_c(-\infty)=1,\ \ \varphi_c(+\infty)=0.
\ee
Furthermore, the function $\varphi_c$ is decreasing in $\R$ and unique up to shifts. Now, for pro\-blem~(\ref{cauchy}) under the sole assumption (\ref{frontlike}), when --instead of (\ref{slow})-- $u_0(x)$ is equivalent as $x\to+\infty$ to a multiple of $e^{-\alpha x}$ with $0<\alpha<\alpha^*=\sqrt{f'(0)}$, then $u(t,x)$ converges to a finite shift of the front $\varphi_c(x-ct)$ as $t\to+\infty$, where
$$c=\alpha+\frac{f'(0)}{\alpha}>c^*.$$
On the other hand, when $u_0(x)=O(e^{-\alpha^*x})$ as $x\to+\infty$, then $u(t,x)$ behaves as $t\to+\infty$ like $\varphi_{c^*}(x+m(t)-c^*t)$, where the phase shift $m(t)$ satisfies $m(t)=O(\ln t)$ (the particular case when $u_0=0$ on $[0,+\infty)$ was first treated in the seminal paper of Kolmogorov, Petrovski and Piskunov \cite{kpp}). In these two situations, the ``location" of the solution $u$ at large time moves at a finite speed, in the sense that, for any $\lambda\in(0,1)$ and any family of real numbers $x_{\lambda}(t)$ such that $u(t,x_{\lambda}(t))=\lambda$, then $x_{\lambda}(t)/t$ converges as $t\to+\infty$ to a positive constant. This constant asymptotic speed is equal to $c=\alpha+f'(0)/\alpha$ in the first case, it is equal to $c^*$ in the second case. We refer to \cite{bhm,b,ev,k,la,l,mk,r,u} for a much more complete picture and more detailed results concerned with exponentially bounded initial conditions~$u_0$ (see also~\cite{hn} for further estimates of possibly different propagation speeds when the initial condition $u_0$ is just assumed to be trapped as $x\to+\infty$ between two exponentially decaying functions with different exponents).\par
Most of these results also hold for more general reaction-diffusion equations with Fisher-KPP nonlinearities: the global stability of traveling fronts in infinite cylinders with shear flows has been shown in \cite{mr}, the convergence to pulsating traveling fronts for reaction-diffusion-advection equations in periodic media has been established in~\cite{hr}. Further estimates about the {\sl finite} spreading speeds associated with compactly supported or exponentially decaying initial conditions have been obtained in \cite{aw} for homogeneous equations in $\R^N$, in~\cite{bmr} for one-dimensional periodic equations, in \cite{bhn2,kr,lyz,nrx,w,x2} for higher-dimensional equations set in periodic media, or in \cite{ckor} for equations with more general advection. We also mention that various definitions and estimates of the spreading speeds have been given in time almost and space periodic media \cite{hush} and in general domains~\cite{bhn,bhn3}. Much work in the mathematical literature has also been devoted to these propagation issues for other types of nonlinearities, like for more general positive non-KPP nonlinearities \cite{aw,ev,hr,r2,w,z2} (in the specific case $f(s)=s^m(1-s)$ with $m>1$, the situation is quite different from the KPP case since traveling fronts with non-minimal speeds have algebraic decay at $+\infty$ and algebraically decaying initial condition may travel with finite or infinite speeds, according to the decay rate, see \cite{ksm,nb,sm}), combustion-type nonlinearities \cite{blr,hs,k1,k2,mnrr,r1,r2,z1} or bistable type nonlinearities \cite{fm,lx,mnt,r2,x1,z1}.\par
Let us now come back to our problem (\ref{cauchy}) with a Fisher-KPP nonlinearity $f$ satisfying~(\ref{fkpp}) and an initial condition $u_0$ decaying slowly to $0$ at $+\infty$ in the sense of (\ref{frontlike}) and~(\ref{slow}). Apart from some estimates on time-exponentially accelerating solutions when $f(s)=s(1-s)$ and when the initial condition $u_0$ converges algebraically to~$0$ and~$1$ at~$\pm\infty$ (see Theorem~2 in \cite{la}, see also \cite{nb} for similar problems set in the half-line $[0,+\infty)$ with nonincreasing algebraically decaying $u_0$), the important issues of the ``location" and of the shape of the profile of $u$ at large time for general KPP functions $f$ and under the general assumptions (\ref{frontlike}) and (\ref{slow}), without any monotonicity or algebraic decay assumption at initial time, have not yet been dealt with, despite of their relevance and simplicity. We shall see that completely different and new propagation phenomena occur: firstly, the level sets of any given value $\lambda\in(0,1)$ (namely, the time-dependent sets of real numbers~$x$ such that $u(t,x)=\lambda$) travel infinitely fast as $t\to+\infty$, are quantitatively estimated in terms of the behavior of $u_0$ at $+\infty$ and can be chosen as large as wanted, and, secondly, the solution~$u$ becomes uniformly flat as~$t\to+\infty$.\par
Our first main result is thus concerned with the large-time behavior of the level sets of the solution $u$, under assumptions (\ref{frontlike}) and (\ref{slow}). Before stating the theorem, we need to introduce a few notations. For any $\lambda\in(0,1)$ and $t\ge 0$, denote by
$$E_{\lambda}(t)=\{x\in\R,\ u(t,x)=\lambda\}$$
the level set of $u$ of value $\lambda$ at time $t$. For any subset $A\subset(0,1]$, we set
$$u_0^{-1}\big\{A\big\}=\{x\in\R,\ u_0(x)\in A\}$$
the inverse image of $A$ by $u_0$.

\begin{theo}\label{th1}
Let $u$ be the solution of $(\ref{cauchy})$, where $f$ satisfies $(\ref{fkpp})$ and the initial condition $u_0:\R\to[0,1]$ satisfies~$(\ref{frontlike})$ and~$(\ref{slow})$.\par
$a)$ Then
$$\lim_{x\to+\infty}u(t,x)=0\hbox{ for all }t\ge 0,\ \hbox{ and }\ \liminf_{x\to-\infty}u(t,x)\to 1\hbox{ as }t\to+\infty.$$\par
$b)$ For any given $\lambda\in(0,1)$, there is a real number $t_{\lambda}\ge 0$ such that $E_{\lambda}(t)$ is compact and non-empty for all $t\ge t_{\lambda}$, and
\be\label{infinite}
\lim_{t\to+\infty}\frac{\min E_{\lambda}(t)}{t}=+\infty.
\ee\par
$c)$ Assume that there exists $\xi_0\in\R$ such that $u_0$ is of class $C^2$ and nonincreasing on $[\xi_0,+\infty)$, and $u''_0(x)=o(u_0(x))$ as $x\to+\infty$.\footnote{Notice that when $u_0$ satisfies (\ref{frontlike}) and $u''_0(x)=o(u_0(x))$ as $x\to+\infty$, then the property $u'_0(x)=o(u_0(x))$ is automatically fulfilled, whence~(\ref{slow}) holds. Indeed, the function $g=u'_0/u_0$ satisfies $g'=u''_0/u_0-g^2$ in $[\xi_0,+\infty)$. Now, if there exist $\epsilon>0$ and an increasing sequence $(y_n)_{n\in\N}$ in $[\xi_0,+\infty)$ such that $\lim_{n\to+\infty}y_n=+\infty$ and $|g(y_n)|\ge\epsilon$ for all $n$, then there exists $N\in\N$ such that $u''_0/u_0\le\epsilon^2/2$ on $[y_N,+\infty)$ and $g'(y_n)\le-\epsilon^2/2$ for all $n\ge N$. Up to extraction of a subsequence, two cases may occur: either $g(y_n)\le-\epsilon$ for all $n$, or $g(y_n)\ge\epsilon$ for all $n$. In the first case, $g$ is decreasing on $[y_N,+\infty)$ and $g'\le-\epsilon^2/2$ on $[y_N,+\infty)$, whence $g(+\infty)=g'(+\infty)=-\infty$. Therefore, $g'(x)\le-g(x)^2/2$ for large $x$ and a contradiction follows by integrating this inequality. In the second case, then, for all $n\ge N$, $g'\le\epsilon^2/2-g^2$ on $[y_N,y_n]$ and $g$ is decreasing and $g\ge\epsilon$ on this interval. A contradiction follows by integrating between $y_N$ and $y_n$ and passing to the limit as $n\to+\infty$.} Then, for any $\lambda\in(0,1)$, $\epsilon\in(0,f'(0))$, $\gamma>0$ and $\Gamma>0$, there exists $T_{\lambda,\epsilon,\gamma,\Gamma}\ge t_{\lambda}$ such that
\be\label{largetime}
\forall\, t\ge T_{\lambda,\epsilon,\gamma,\Gamma},\quad
E_{\lambda}(t)\subset u_0^{-1}\Big\{[\gamma\,e^{-(f'(0)+\epsilon)t},\Gamma\,e^{-(f'(0)-\epsilon)t}]\Big\}.
\ee
\end{theo}

Part~b) of Theorem~\ref{th1} simply says that the level sets $E_{\lambda}(t)$ of all level values $\lambda\in(0,1)$ move infinitely fast as $t\to+\infty$, in the average sense (\ref{infinite}). As already announced above, this property is in big contrast with the finiteness of the propagation speeds of solutions which are exponentially bounded as $x\to+\infty$ at initial time. Moreover, Theorem~\ref{th1} actually yields the following corollary, which states that the level sets $E_{\lambda}(t)$ can be located as far to the right as wanted, provided that the initial condition is well chosen.

\begin{cor}\label{cor1}
Under the assumptions and notations of Theorem~$\ref{th1}$, the following holds: given any function $\xi:[0,+\infty)\to\R$ which is locally bounded, there are initial conditions~$u_0$ such that, for any given $\lambda\in(0,1)$,
$$\min E_{\lambda}(t)\ge\xi(t)\ \hbox{ for all }t\hbox{ large enough}.$$
\end{cor}

Let us now comment on the quantitative estimates given in part~c) of Theorem~\ref{th1}. Observe first that the real numbers $\gamma\,e^{-(f'(0)+\epsilon)t}$ and $\Gamma\,e^{-(f'(0)-\epsilon)t}$ belong to $(0,\sup_{\R}u_0)$ for large $t$, giving a meaning to~(\ref{largetime}). Now, if there is $a\in\R$ such that $u_0$ is {\sl strictly} decreasing on $[a,+\infty)$, then the inclusions~(\ref{largetime}) mean that
\be\label{largetimebis}
u_0^{-1}\big(\Gamma\,e^{-(f'(0)-\epsilon)t}\big)\le\min E_{\lambda}(t)\le\max E_{\lambda}(t)\le u_0^{-1}\big(\gamma\,e^{-(f'(0)+\epsilon)t}\big).
\ee
for large $t$, where, here, $u_0^{-1}:(0,u_0(a)]\to[a,+\infty)$ denotes the reciprocal of the restriction of the function $u_0$ on the interval $[a,+\infty)$. Furthermore, if $u_0$ is nonincreasing over the whole real line $\R$, then the derivative $u_x(t,x)$ is negative for all $t>0$ and $x\in\R$, from the strong parabolic maximum principle. Therefore, $E_{\lambda}(t)$ is either empty or a singleton as soon as $t>0$. In particular, $E_{\lambda}(t)=:\{x_{\lambda}(t)\}$ for all $t>t_{\lambda}$, and the maps $t\mapsto x_{\lambda}(t)$ are of class $C^1(t_{\lambda},+\infty)$ from the implicit function theorem. The inequalities (\ref{largetimebis}) then provide lower and upper bounds of $x_{\lambda}(t)$ for large $t$. However, it is worth pointing out that Theorem~\ref{th1} is valid for general initial conditions which decay slowly to $0$ at $+\infty$ but which may not be globally nonincreasing.\par
Notice that, given any two values $\lambda$ and $\lambda'$ in $(0,1)$, assertion (\ref{largetime}) implies that the level sets $E_{\lambda}(t)$ and $E_{\lambda'}(t)$ are both included in the same family of moving sets $u_0^{-1}\Big\{[\gamma\,e^{-(f'(0)+\epsilon)t},\Gamma\,e^{-(f'(0)-\epsilon)t}]\Big\}$ for $t$ large enough. Although the dependence on $\lambda$ and $\lambda'$ is not explicit, the estimates (\ref{largetime}) are still sufficient to derive explicit large-time asymptotic equivalents (or logarithmic-type asymptotic equivalents) of the positions of all level sets~$E_{\lambda}(t)$ for the main aforementioned examples (see the end of this section).

\begin{rem}{\rm The estimate (\ref{infinite}) and the ``lower bound" of  $\min E_{\lambda}(t)$ in (\ref{largetime}), that is
\be\label{largetimelower}
E_{\lambda}(t)\subset u_0^{-1}\Big\{(0,\Gamma\,e^{-(f'(0)-\epsilon)t}]\Big\}
\ee
for $t$ large enough hold under weaker assumptions on $f$: namely, for~(\ref{infinite}) and (\ref{largetimelower}) to hold, $f$ can simply be assumed to vanish at~$0$ and~$1$, to be positive on $(0,1)$ and to have a positive derivative $f'(0)>0$ at $0$. We refer to the proof of Theorem~\ref{th1} in Section~\ref{sec3} for the details. Similarly, Corollary~\ref{cor1} and the lower bound in (\ref{largetimebis}) if $u_0$ is decreasing in $[a,+\infty)$ also hold under these weaker assumptions.}
\end{rem}

Next, let us give a further interpretation of Theorem~\ref{th1} and an insight about the main underlying ideas of the proofs. Let $\lambda\in(0,1)$ and let $t\mapsto x_{\lambda}(t)$ be any map such that $x_{\lambda}(t)\in E_{\lambda}(t)$ for large $t$, that is $u(t,x_{\lambda}(t))=\lambda$. Formula (\ref{largetime}) with $\gamma=\Gamma=\lambda$ can be rewritten as
\be\label{xlambdaeps}
u_0(x_{\lambda}(t))\,e^{(f'(0)-\epsilon)t}\le\lambda\le u_0(x_{\lambda}(t))\,e^{(f'(0)+\epsilon)t}
\ee
for $t$ large enough. Roughly speaking, this means that the real numbers $x_{\lambda}(t)$ are then asymptotically given, in the above approximate sense, by the solution of the family of decoupled ODE's
\be\label{eqU}\left\{\baa{rcl}
\displaystyle{\frac{{\rm{d}}U(t;x)}{{\rm{d}}t}} & = & f'(0)U(t;x),\vspace{5pt}\\
U(0;x) & = & u_0(x),\eaa\right.
\ee
parameterized by $x\in\R$, and then, say, by solving $U(t;x_{\lambda}(t))=\lambda$. In other words, the behavior of $u(t,x)$ at large time is dominated by the reaction term, that is to say that the diffusion term plays in some sense a negligible role as compared to the growth by reaction. Our method of proof in Section~\ref{sec3} is based on the rigorous formulation of these observations.\par
Actually, (\ref{eqU}) is a good approximation of (\ref{largetime}) and (\ref{xlambdaeps}) up to the $\pm\epsilon$ terms in the exponents of the exponential factors. It turns out that, under additional assumptions, the~$\pm\epsilon$ terms can be dropped and more precise estimates than (\ref{largetime}) can be established, as the following theorem shows:

\begin{theo}\label{th1bis}
Let the general assumptions of Theorem~$\ref{th1}$ be fulfilled. If, on the one hand, one assumes moreover that $u_0$ is convex for large $x$, or that $u_0''(x)=O(u_0(x)^{1+\beta})$ as $x\to+\infty$ for some $\beta>0$, then, for all $\lambda\in(0,1)$, there exist $\Gamma_{\lambda}>0$ and a time $\underline{T}_{\lambda}\ge t_{\lambda}$ such that
$$\forall\,t\ge\underline{T}_{\lambda},\quad E_{\lambda}(t)\subset u_0^{-1}\Big\{(0,\Gamma_{\lambda}\,e^{-f'(0)t}]\Big\},$$
where the positive real numbers $\Gamma_{\lambda}$ can be chosen independently of $\lambda$ when $\lambda>0$ is small. On the other hand, if  there exist $\mu>0$ and $\beta\ge\nu>0$ such that
$$f(s)\le f'(0)s-\mu\,s^{1+\nu}\ \hbox{ for all }s\in[0,1],\footnote{Notice that this condition is satisfied with $\nu=1$ for instance if $f$ is of class $C^2([0,1])$ and strictly concave, that is $f''(s)<0$ for all $s\in[0,1]$.}$$
and $u_0''(x)=O(u_0(x)^{1+\beta})$ as $x\to+\infty$, then there exists a constant $c>0$ such that, for all $\lambda\in(0,1)$, there exists a time $\overline{T}_{\lambda}\ge t_{\lambda}$ such that
\be\label{upperprecise}
\forall\,t\ge\overline{T}_{\lambda},\quad E_{\lambda}(t)\subset u_0^{-1}\Big\{[c\,\lambda\,e^{-f'(0)t},1]\Big\}.
\ee
\end{theo}

Our last result provides additional information on the global aspect of the graph of the functions $u(t,\cdot)$ at large time, under some appropriate assumptions.

\begin{theo}\label{th2}
Assume that $f$ satisfies $(\ref{fkpp})$ and $f'(s)\le f(s)/s$ for all $s\in(0,1]$. Assume that $u_0:\R\to[0,1]$ satisfies~$(\ref{frontlike})$ and~$(\ref{slow})$, and that $u'_0/u_0\in L^p(\R)\cap C^{2,\theta}(\R)$ for some $p\in(1,+\infty)$ and $\theta\in(0,1)$. Then the solution $u$ of $(\ref{cauchy})$ is such that
\be\label{flat}
\|u_x(t,\cdot)\|_{L^{\infty}(\R)}\to 0\hbox{ as }t\to+\infty.
\ee
\end{theo}

In other words, the solution $u$ becomes uniformly flat at large time, in the sense that there are no regions on which the solution has spatial gradients~$u_x$ which are bounded away from $0$ as $t\to+\infty$. The extra-assumption on $f$ means that the growth rate $f(s)/s$ is nonincreasing with respect to $s$ (this is the case for instance if $f$ is concave). K.~Uchiyama proved in \cite{u} (see Theorem~8.5 in \cite{u}) that $u(t,x+\max E_{1/2}(t))\to1/2$ locally uniformly in $x$ as $t\to+\infty$ (the value $1/2$ could clearly be replaced by any real number $0<\lambda<1$). Standard parabolic estimates then imply that $u_x(t,x+\max E_{\lambda}(t))\to 0$ locally uniformly in $x$ as $t\to+\infty$. Formula (\ref{flat}) yields more results since it implies that the convergence of $u_x$ to $0$ as $t\to+\infty$ takes place not only around the position $\max E_{\lambda}(t)$, but also around any point of $E_{\lambda}(t)$ and uniformly both with respect to all these points and with respect to all values $\lambda\in(0,1)$. Notice that since $u_0$ is not assumed to be nonincreasing, the sets $E_{\lambda}(t)$ may then not reduce to singletons and the quantities $\max E_{\lambda}(t)-\min E_{\lambda}(t)$ may not be bounded, all the more as the profile $u(t,\cdot)$ does become flater and flater as $t\to+\infty$. However, we use in Theorem~\ref{th2} an additional smoothness assumption on $u_0$ and an integrability condition, that is $u'_0/u_0\in L^p(\R)$ for some $p\in(1,+\infty)$.\footnote{With~(\ref{frontlike}), the assumption $u'_0/u_0\in L^p(\R)$ is not so restrictive since it means that $u'_0\in L^p(-\infty,0)$ and $u'_0/u_0\in L^p(0,+\infty)$ (this last condition is satisfied for instance for all four examples~(\ref{ex1}-\ref{ex4}) listed above). Since $u'_0/u_0$ is continuous and converges to $0$ at $\pm\infty$, it then follows that $u'_0/u_0\in L^q(\R)$ for all $q\in[p,+\infty]$.} This assumption actually allows us to prove more than (\ref{flat}), since we show in Section~\ref{sec3} that
\be\label{flatbis}
\left\|\frac{u_x}{u}(t,\cdot)\right\|_{L^{\infty}(\R)}\to 0\ \hbox{ as }t\to+\infty.
\ee
This property yields large-time exponential estimates which are uniform with respect to any origin $y$: for all $\epsilon>0$, there exists a time $T_{\epsilon}\ge0$ such that, for all $t\ge T_{\epsilon}$ and for all $y\le x\in\R$,
$$u(t,x)\ge u(t,y)\,e^{-\epsilon(x-y)}.$$

\begin{rem}{\rm Assume here that $u_0$ is nonincreasing and set, for each $t\ge 0$, $p^+(t)=0$ and $p^-(t)=u(t,-\infty)$ (which exists since $u(t,\cdot):\R\to[0,1]$ is nonincreasing). At each time $t$, the function $u(t,\cdot)$ connects $p^-(t)$ to $p^+(t)$ in the sense that $u(t,\pm\infty)=p^{\pm}(t)$. However, an important consequence of Theorem~\ref{th2} is that, under the assumptions of this theorem, the solution $u$ is {\sl not} a generalized transition wave between $p^-$ and $p^+$ in the sense of \cite{bh4,bh5}. Namely, it follows from~(\ref{flat}) that there is no real-valued map $t\mapsto x(t)$ such that $u(t,x)-p^{\pm}(t)\to 0$ as $x-x(t)\to\pm\infty$, uniformly with respect to $t\ge 0$. Indeed, if there were such real numbers $x(t)$, then $u$ would have a ``width" which should be bounded uniformly with respect to $t$. But this is ruled out from (\ref{flat}) and part~a) of Theorem~\ref{th1}.}
\end{rem}

To complete this section, let us now apply the estimates of Theorem~\ref{th1} to the four main examples given in (\ref{ex1}-\ref{ex4}). If $u_0$ satisfies (\ref{frontlike}) and
$$u_0(x)\sim C\,e^{-\alpha\,x/\ln x}\hbox{ as }x\to+\infty$$
with $\alpha,\ C>0$, and if $u''_0(x)=o(u_0(x))$ as $x\to+\infty$ (this condition is fulfilled if $u_0(x)=C\,e^{-\alpha\,x/\ln x}$ for large $x$), then it follows from (\ref{largetime}) that the positions of the level sets $E_{\lambda}(t)$ of any level value $\lambda\in(0,1)$ satisfy
$$\min E_{\lambda}(t)\sim \max E_{\lambda}(t)\sim f'(0)\,\alpha^{-1}\,t\,\ln t\ \hbox{ as }t\to+\infty.$$
Notice that the equivalent of the $\min$ and $\max$ of $E_{\lambda}(t)$ does not depend on $\lambda$. However, under the assumptions of Theorem~\ref{th2} and if one assumes that $u_0$ is nonincreasing, there holds $\lim_{t\to+\infty}|x_{\lambda}(t)-x_{\lambda'}(t)|=+\infty$ for any $\lambda\neq\lambda'\in(0,1)$, where, in this case, $E_{\lambda}(t)=\{x_{\lambda}(t)\}$ and $E_{\lambda'}(t)=\{x_{\lambda'}(t)\}$ for large $t$.\par
Now, if $u_0$ satisfies (\ref{frontlike}) and
$$u_0(x)\sim C\,e^{-\beta\,x^{\alpha}}\hbox{ as }x\to+\infty$$
with $\alpha\in(0,1)$ and $\beta,\ C>0$, and if $u''_0(x)=o(u_0(x))$ as $x\to+\infty$ (this condition is fulfilled if $u_0(x)=C\,e^{-\beta\,x^{\alpha}}$ for large $x$), then it follows from (\ref{largetime}) that the positions of the level sets $E_{\lambda}(t)$ of any level value $\lambda\in(0,1)$ are asymptotically algebraic and superlinear as $t\to+\infty$, in the sense that
$$\min E_{\lambda}(t)\sim \max E_{\lambda}(t)\sim f'(0)^{1/\alpha}\,\beta^{-1/\alpha}\,t^{1/\alpha}\ \hbox{ as }t\to+\infty.$$\par
If $u_0$ satisfies (\ref{frontlike}) and
\be\label{ex3bis}
u_0(x)\sim C\,x^{-\alpha}\hbox{ as }x\to+\infty
\ee
with $\alpha,\ C>0$, and if $u''_0(x)=o(u_0(x))$ as $x\to+\infty$ (this condition is fulfilled if $u_0(x)=C\,x^{-\alpha}$ for large $x$), then it follows from (\ref{largetime}) that the positions of the level sets $E_{\lambda}(t)$ of any level value $\lambda\in(0,1)$ move exponentially fast as $t\to+\infty$:
\be\label{logequiv}
\ln(\min E_{\lambda}(t))\sim \ln(\max E_{\lambda}(t))\sim f'(0)\,\alpha^{-1}\,t\hbox{ as }t\to+\infty.
\ee
We mention that X.~Cabr\'e and J.-M.~Roquejoffre \cite{cr} just established similar estimates for the level sets of the solutions $u$ of equations of the type $u_t=Au+f(u)$, where $f$ is concave,~$u_0$ is compactly supported or monotone one-sided compactly supported, and the operator~$A$ is the generator of a Feller semi-group (a typical example is the fractional Laplacian $A=(-\Delta)^{\rho}$ with $0<\rho<1$). In this case, the asymptotic exponential spreading of the level sets follows from the algebraic decay of the kernels associated to the operators~$A$. For problem (\ref{cauchy}) and (\ref{ex3bis}), if $u_0$ and $f$ satisfy the additional conditions of Theorem~\ref{th1bis}, that is if $f(s)\le f'(0)s-\mu\,s^{1+\nu}$ in $[0,1]$ and $u_0''(x)=O(u_0(x)^{1+\beta})$ as $x\to+\infty$ for some $\mu>0$ and $\beta\ge\nu>0$,\footnote{This is the case with $\beta=\nu=1$ and $0<\alpha\le2$ if $f$ is $C^2$ and strictly concave on $[0,1]$ and if $u_0(x)=C\,x^{-\alpha}$ for large $x$.} then Theorem~\ref{th1bis} says that one can be more precise than (\ref{logequiv}), in the sense that, for every $\lambda\in(0,1)$,
$$C_{\lambda}\,e^{f'(0)\,\alpha^{-1}\,t}\le\min E_{\lambda}(t)\le\max E_{\lambda}(t)\le C'_{\lambda}\,e^{f'(0)\,\alpha^{-1}\,t}$$
for some constants $0<C_{\lambda}\le C'_{\lambda}$ and for $t$ large enough.\par
Lastly, if $u_0$ satisfies (\ref{frontlike}) and
$$u_0(x)\sim C\,(\ln x)^{-\alpha}\hbox{ as }x\to+\infty$$
with $\alpha,\ C>0$, and if $u''_0(x)=o(u_0(x))$ as $x\to+\infty$ (this condition is fulfilled if $u_0(x)=C\,(\ln x)^{-\alpha}$ for large $x$), then it follows from (\ref{largetime}) that the positions of the level sets $E_{\lambda}(t)$ of any level value $\lambda\in(0,1)$ move doubly-exponentially fast as $t\to+\infty$:
$$\ln(\ln(\min E_{\lambda}(t)))\sim \ln(\ln(\max E_{\lambda}(t)))\sim f'(0)\,\alpha^{-1}\,t\hbox{ as }t\to+\infty.$$
Obviously, it is immediate to see that these examples can be generalized. For instance,~$\min E_{\lambda}(t)$ and $\max E_{\lambda}(t)$ move $m$-exponentially fast with $m\ge 2$, that is $\ln^{\circ(m)}(\min E_{\lambda}(t))$ and $\ln^{\circ(m)}(\max E_{\lambda}(t))$ behave linearly in $t$ as $t\to+\infty$ if $u_0$ behaves like a negative power of $\ln^{\circ(m-1)}(x)$ as $x\to+\infty$, where $\ln^{\circ(1)}=\ln$ and $\ln^{\circ(k)}=\ln\circ\ln^{\circ(k-1)}$ for $k\ge 2$. More precise estimates also hold under the assumptions of Theorem~\ref{th1bis}.
\hfill\break

\noindent{\bf{Outline of the paper.}} The remaining part of the article is organized as follows: Section~\ref{sec2} is devoted to the proof of the large-time asymptotics of the positions of the level sets $E_{\lambda}(t)$ of $u$, that is Theorem~\ref{th1}, Corollary~\ref{cor1} and Theorem~\ref{th1bis}. In Section~\ref{sec3}, we prove Theorem~\ref{th2}, that is the solutions $u$ become uniformly flat at large times.


\SE{Motion of the level sets}\label{sec2}

This section is concerned with the proof of Theorem~\ref{th1} and its consequences: Corollary~\ref{cor1} and Theorem~\ref{th1bis}. Parts~a) and~b)  of Theorem~\ref{th1} follow from elementary comparisons with traveling fronts. Part~c) is more involved and is based on the comparison of the solution $u$ with the solutions of systems of ODEs which approximate~(\ref{eqU}) and which are sub- and supersolutions for the original problem (\ref{cauchy}). The proof of Theorem~\ref{th1bis} is based on more involved estimates and refinments of the comparison functions which are used in the proof of part~c) of Theorem~\ref{th1}. One can then drop the $\pm\epsilon$ terms in (\ref{largetime}).\hfill\break

\noindent{\bf{Proof of Theorem~\ref{th1}.}} Parts~a) and~b) follow immediately from the classical results on the stability of traveling fronts with finite speed. We sketch here the main ideas for the sake of completeness. First, for any sequence $(x_n)_{n\in\N}$ such that $\lim_{n\to+\infty}x_n=+\infty$, the functions $u_n(t,x)=u(t,x+x_n)$ converge locally uniformly in $[0,+\infty)\times\R$, up to extraction of a subsequence, to a classical solution $u_{\infty}$ of (\ref{cauchy}) with initial condition $u_{\infty}(0,\cdot)=0$. Therefore, $u_{\infty}=0$ in $[0,+\infty)\times\R$. The uniqueness of the limit implies that $u(t,x)\to0$ as $x\to+\infty$ for each $t\ge 0$.\par
Now, from (\ref{frontlike}) and (\ref{slow}), for each $c>2\sqrt{f'(0)}$, there exist a uniformly continuous decreasing function $v_0:\R\to[0,1]$ and $A\in\R$ such that $0<v_0\le u_0$ in $\R$ and $v_0(x)=e^{-\alpha_cx}$ for all $x\ge A$, where $\alpha_c=(c-\sqrt{c^2-4f'(0)})/2>0$. Denote by $v$ the solution of the Cauchy problem (\ref{cauchy}) with initial condition $v_0$. It is well-known (see Section~\ref{intro}) that
$$\|v(t,\cdot)-\varphi_c(\cdot-ct+\xi)\|_{L^{\infty}(\R)}\to 0\ \hbox{ as }t\to+\infty$$
for some $\xi\in\R$, where $\varphi_c$ solves (\ref{varphic}). But the maximum principle yields $v(t,x)\le u(t,x)\ (\le 1)$ for all $t\ge 0$ and $x\in\R$. Therefore,
$$\liminf_{x\to-\infty}u(t,x)\to 1\ \hbox{ as }\ t\to+\infty.$$
It then follows, that for any $\lambda\in(0,1)$, there exists a time $t_{\lambda}\ge 0$ such that
$$\liminf_{x\to-\infty}u(t,x)>\lambda>0=u(t,+\infty)$$
for all $t\ge t_{\lambda}$. Since all functions $x\mapsto u(t,x)$ are continuous, one concludes that $E_{\lambda}(t)$ is a non-empty compact set for all $t\ge t_{\lambda}$. Furthermore, since $c$ can be chosen arbitrarily large, one concludes that $\min E_{\lambda}(t)/t\to+\infty$ as $t\to+\infty$.

\begin{rem}{\rm The conclusions of parts~a) and b) still hold when, instead of (\ref{fkpp}), the function $f$ is only assumed to vanish at $0$ and $1$, to be positive on $(0,1)$ and to have a positive derivative~$f'(0)>0$ at~$0$. Indeed, as far as the lower estimates are concerned, it is sufficient to replace $f$ by any $C^1$ function $\widetilde{f}:[0,1]\to\R$ satisfying (\ref{fkpp}) and such that $\widetilde{f}\le f$ in $[0,1]$.}
\end{rem}

Let us  now turn to the proof of part~c) of Theorem~\ref{th1}. We begin with the upper bounds of $\max E_{\lambda}(t)$ for large $t$. We here fix $\lambda\in(0,1)$, $\epsilon>0$ and $\gamma>0$. We shall prove that
\be\label{claim1}
E_{\lambda}(t)\subset u_0^{-1}\Big\{[\gamma\,e^{-(f'(0)+\epsilon)t},1]\Big\}
\ee
for $t$ large enough. To do so, first choose $\xi_1\in[\xi_0,+\infty)$ so that $|u''_0(x)|\le(\epsilon/2)\,u_0(x)$ for all $x\ge\xi_1$ (remember that $u_0$ is of class $C^2$ on $[\xi_0,+\infty)$ and that $u_0''(x)=o(u_0(x))$ as $x\to+\infty$). Set $\rho=f'(0)+\epsilon/2$ and, for all $(t,x)\in[0,+\infty)\times[\xi_1,+\infty)$,
\be\label{defoveru}
\overline{u}(t,x)=\min\left(\frac{u_0(x)\,e^{\rho t}}{u_0(\xi_1)},1\right).
\ee
Observe that $u_0(x)\le\overline{u}(0,x)$ for all $x\in[\xi_1,+\infty)$, and that $u(t,\xi_1)\le 1=\overline{u}(t,\xi_1)$ for all $t\in[0,+\infty)$. Let us now check that $\overline{u}$ is a supersolution of the equation satisfied by $u$, in the set $[0,+\infty)\times[\xi_1,+\infty)$. Since $u\le 1$ and $f(1)=0$, it is sufficient to check it when $\overline{u}<1$. If $(t,x)\in[0,+\infty)\times[\xi_1,+\infty)$ and $\overline{u}(t,x)<1$, then it follows from (\ref{fkpp}) that
$$\baa{rcl}
\overline{u}_t(t,x)-\overline{u}_{xx}(t,x)-f(\overline{u}(t,x)) & \ge & \displaystyle{\frac{1}{u_0(\xi_1)}}\times\big(\rho\,u_0(x)-u_0''(x)-f'(0)\,u_0(x)\big)e^{\rho t}\vspace{5pt}\\
& = & \displaystyle{\frac{1}{u_0(\xi_1)}}\times\left(\displaystyle{\frac{\epsilon}{2}}\,u_0(x)-u_0''(x)\right)e^{\rho t}\vspace{5pt}\\
& \ge & 0\eaa$$
due to the choice of $\xi_1$. The parabolic maximum principle then implies that
\be\label{ineq1}
\forall\,(t,x)\in[0,+\infty)\times[\xi_1,+\infty),\quad u(t,x)\le\overline{u}(t,x)\le u_0(\xi_1)^{-1}u_0(x)\,e^{\rho t}.
\ee
For all $t\ge t_{\lambda}$ (so that $E_{\lambda}(t)$ is not empty) and for all $y\in E_{\lambda}(t)$, there holds
\be\label{ineq2}
\Big(\hbox{either }y<\xi_1\Big)\hbox{ or }\Big(y\ge\xi_1\hbox{ and }\lambda=u(t,y)\le u_0(\xi_1)^{-1}u_0(y)\,e^{\rho t}\Big).
\ee
In all cases, one gets that
\be\label{ineq3}
\forall\, t\ge t_{\lambda},\ \forall\,y\in E_{\lambda}(t),\quad u_0(y)\ge\min\left(\eta,\lambda\,u_0(\xi_1)\,e^{-\rho t}\right),
\ee
where $\eta=\inf_{(-\infty,\xi_1)}u_0>0$. Since $\rho=f'(0)+\epsilon/2$, there exists then a time $\overline{T}_{\lambda,\epsilon,\gamma}\ge t_{\lambda}$ such that
\be\label{claim1bis}
\forall\, t\ge\overline{T}_{\lambda,\epsilon,\gamma},\ \forall\,y\in E_{\lambda}(t),\quad u_0(y)\ge\gamma\,e^{-(f'(0)+\epsilon)t},
\ee
which gives (\ref{claim1}).\hfill\break\par

We now turn to the proof of the lower bounds of $\min E_{\lambda}(t)$ for large $t$. We fix $\lambda\in(0,1)$, $\epsilon\in(0,f'(0))$ and $\Gamma>0$. We claim that
\be\label{claim2}
E_{\lambda}(t)\subset u_0^{-1}\Big\{(0,\Gamma\,e^{-(f'(0)-\epsilon)t}]\Big\}
\ee
for $t$ large enough. Here, we do not need to assume (\ref{fkpp}), but, instead, we simply assume that $f(0)=f(1)=0$, $f(s)>0$ on $(0,1)$ and $f'(0)>0$. However, we recall that there exist $\delta>0$, $s_0\in(0,1)$ and $M\ge 0$ such that $f(s)\ge f'(0)\,s-M\,s^{1+\delta}$ for all $s\in[0,s_0]$.\par
To do the proof of (\ref{claim2}), let us first introduce a few notations. Choose $\rho\in\R$ so that
$$f'(0)-\epsilon<\rho<f'(0)\ \hbox{ and }\ \rho\,(1+\delta)>f'(0).$$
Next, choose $\xi_2\in[\xi_0,+\infty)$ such that
\be\label{xi3}
\forall\,x\ge\xi_2,\quad|u''_0(x)|\le\min\left(f'(0)-\rho,\frac{\rho\,(1+\delta)-f'(0)}{2\,(1+\delta)}\right)\times u_0(x).
\ee
Then, define
\be\label{kappa}
\kappa=\inf_{(-\infty,\xi_2)}u_0>0,\ s_1=\min(s_0,\kappa)>0\hbox{ and }B=\max\left(s_1^{-\delta},\frac{2\,M}{\rho\,(1+\delta)-f'(0)}\right)\!\!>0.
\ee
Let $g$ be the function defined in $\R_+$ by
$$g(s)=s-B\,s^{1+\delta}.$$
Observe that
$$g(s)\le 0\hbox{ for all }s\ge s_1\hbox{ and }g(s)\le s_1\hbox{ for all }s\ge 0.$$
Lastly, denote by $\underline{u}$ the function defined by
\be\label{defunderu}
\underline{u}(t,x)=\max\left(g(u_0(x)\,e^{\rho t}),0\right)=\max\left(u_0(x)\,e^{\rho t}-B\,u_0(x)^{1+\delta}\,e^{\rho(1+\delta)t},0\right)
\ee
for all $(t,x)\in[0,+\infty)\times\R$.\par
Let us check that this function $\underline{u}$ is a subsolution for the Cauchy problem~(\ref{cauchy}). First, there holds $\underline{u}(0,x)\le u_0(x)$ for all $x\in\R$. Since $u\ge 0$ and $f(0)=0$, it is then sufficient to check that $\underline{u}$ is a subsolution of the equation satisfied by~$u$, in the region where $\underline{u}>0$. Let $(t,x)$ be any point in $[0,+\infty)\times\R$ such that $\underline{u}(t,x)>0$. Since $g\le 0$ on $[s_1,+\infty)$, it follows that $u_0(x)\,e^{\rho t}<s_1$, whence $u_0(x)<s_1$ and $x\ge\xi_2$ from~(\ref{kappa}). Furthermore,
\be\label{s10}
0<\underline{u}(t,x)<u_0(x)\,e^{\rho t}<s_1\le s_0<1.
\ee
Thus,
\be\label{lowerf}
f(\underline{u}(t,x))\ge f'(0)\left(u_0(x)\,e^{\rho t}-B\,u_0(x)^{1+\delta}\,e^{\rho(1+\delta)t}\right)-M\,u_0(x)^{1+\delta}\,e^{\rho(1+\delta)t}.
\ee
It follows that
\be\label{subsolution}\baa{l}
\underline{u}_t(t,x)-\underline{u}_{xx}(t,x)-f(\underline{u}(t,x))\vspace{5pt}\\
\quad\le\rho\,u_0(x)\,e^{\rho t}-B\,\rho\,(1+\delta)\,u_0(x)^{1+\delta}\,e^{\rho(1+\delta)t}\vspace{5pt}\\
\quad\quad -u_0''(x) \,e^{\rho t}-B\,(1+\delta)\left(u_0''(x)\,u_0(x)^{\delta}+\delta\,u_0'(x)^2\,u_0(x)^{\delta-1}\right)e^{\rho(1+\delta)t}\vspace{5pt}\\
\quad\quad-f'(0)\left(u_0(x)\,e^{\rho t}-B\,u_0(x)^{1+\delta}\,e^{\rho(1+\delta)t}\right)+M\,u_0(x)^{1+\delta}\,e^{\rho(1+\delta)t}\vspace{5pt}\\
\quad\le\left(\rho\,u_0(x)-u_0''(x)-f'(0)\, u_0(x)\right)e^{\rho t}\vspace{5pt}\\
\quad\quad+\!\left\{B\!\left[\left(f'(0)\!-\!\rho(1+\delta)\right)u_0(x)^{1+\delta}\!-\!(1+\delta)u_0''(x)u_0(x)^{\delta}\right]\!\!+\!Mu_0(x)^{1+\delta}\right\}\!e^{\rho(1+\delta)t}\vspace{5pt}\\
\quad\le0\eaa
\ee
from (\ref{xi3}) and (\ref{kappa}). As a consequence, the maximum principle yields
\be\label{borneinf}
\underline{u}(t,x)\le u(t,x)\ \hbox{ for all }t\ge 0\hbox{ and }x\in\R.
\ee\par
Fix now any real number $\omega$ small enough so that
$$0<\omega<B^{-1/\delta}.$$
This real number $\omega$ does not depend on $\lambda$ nor $\Gamma$, but it depends on $\epsilon$, as well as on~$u_0$ and~$f$. Remember that $t_{\omega}\ge 0$ is such that $E_{\omega}(t)$ is a non-empty compact set for all $t\ge t_{\omega}$. Since $u_0$ is continuous and satisfies (\ref{frontlike}), there exists then a time $\overline{t}_{\omega}\ge t_{\omega}$ such that, for all $t\ge\overline{t}_{\omega}$, the closed set
$$F_{\omega}(t)=\{y\in\R,\ u_0(y)\,e^{\rho t}=\omega\}$$
is non-empty and satisfies $F_{\omega}(t)\subset[\xi_2,+\infty)$. For all $t\ge\overline{t}_{\omega}$, denote
$$y_{\omega}(t)=\min F_{\omega}(t).$$
Since $0<\omega<B^{-1/\delta}\le s_1\le\kappa=\inf_{(-\infty,\xi_2)}u_0$, the function $y_{\omega}:[\overline{t}_{\omega},+\infty)\to[\xi_2,+\infty)$ is nondecreasing and left-continuous, that is $y_{\omega}(t)=y_{\omega}(t^-)=\lim_{s\to t,\,s<t}y_{\omega}(s)$ for all $t>\overline{t}_{\omega}$. Furthermore, since $u_0$ is nonincreasing on $[\xi_2,+\infty)\ (\subset[\xi_0,+\infty))$, the function $y_{\omega}$ is discontinuous at all (and only all) points $t\ge \overline{t}_{\omega}$ for which there exist $a<b\in[\xi_2,+\infty)$ such that $u_0=\omega\,e^{-\rho t}$ on $[a,b]$; if $[a,b]$ denotes the largest such interval, then $a=y_{\lambda}(t)$ and $b=y_{\lambda}(t^+)=\lim_{s\to t,\,s>t}y_{\omega}(s)$.\par
Now, let $\Omega$ be the open set defined by
$$\Omega=\{(t,x)\in(\overline{t}_{\omega},+\infty)\times\R,\ x< y_{\omega}(t)\}.$$
We claim that $\inf_{\overline{\Omega}}u>0$. Observe first that $\partial\Omega$ consists of two parts: the set $\{\overline{t}_{\omega}\}\times(-\infty,y_{\omega}(\overline{t}_{\omega}^+)]$, and the set of all points $(t,x)$ for which $t>\overline{t}_{\omega}$ and $x\in[y_{\omega}(t),y_{\omega}(t^+)]$. If $t>\overline{t}_{\omega}$ and $x\in[y_{\omega}(t),y_{\omega}(t^+)]$, there holds $u_0(x)\,e^{\rho t}=\omega$, whence
$$u(t,x)\ge\omega-B\,\omega^{1+\delta}=g(\omega)>0$$
from~(\ref{defunderu}), (\ref{borneinf}) and the choice of $\omega$. On the other hand, at the time $\overline{t}_{\omega}$, the function $u(\overline{t}_{\omega},\cdot)$ is continuous, positive, and $\liminf_{x\to-\infty}u(\overline{t}_{\omega},x)>0$.\footnote{Indeed, there exist $\eta>0$ and $\xi\in\R$ such that $u_0\ge U_0$ in $\R$, where $U_0=\eta$ on $(-\infty,\xi-1]$, $U_0(x)=\eta\times(\xi-x)$ for all $x\in(\xi-1,\xi)$ and $U_0=0$ on $[\xi,+\infty)$; but the solution $U$ of (\ref{cauchy}) with initial condition $U_0$ satisfies $U\le u$ on $[0,+\infty)\times\R$ and $U(t,-\infty)=\zeta(t)$, where $\zeta$ satisfies $\dot\zeta(t)=f(\zeta(t))$ and $\zeta(0)=\eta$; since $\zeta(t)>0$ for all $t\ge 0$, one concludes that $\liminf_{x\to-\infty}u(t,x)>0$ for all $t\ge 0$.} Thus,
$$\inf_{x\in(-\infty,y_{\omega}(\overline{t}_{\omega}^+)]} u(\overline{t}_{\omega},x)>0.$$
Eventually, there exists $\theta\in(0,1)$ such that $u\ge\theta$ on $\partial\Omega$. Since $f(\theta)>0$, the parabolic maximum principle then implies that
\be\label{Omegatheta}
u(t,x)\ge\theta\ \hbox{ for all }\ (t,x)\in\overline{\Omega}.
\ee\par
Thus, if $\lambda\in(0,\theta)$ and if $x\in E_{\lambda}(t)$ for $t\ge\max(t_{\lambda},\overline{t}_{\omega})$, then
$$x>y_{\omega}(t^+)\ge y_{\omega}(t)\ge\xi_2\ge\xi_0.$$
Since $\Gamma>0$ is fixed and $\rho>f'(0)-\epsilon$, there exists then a time $\underline{T}_{\lambda,\epsilon,\Gamma}\ge\max(t_{\lambda},\overline{t}_{\omega})$ such that,
\be\label{claim2bis}
\forall\, t\ge\underline{T}_{\lambda,\epsilon,\Gamma},\ \forall\,x\in E_{\lambda}(t),\quad u_0(x)\le u_0(y_{\omega}(t))=\omega\,e^{-\rho t}\le\Gamma\,e^{-(f'(0)-\epsilon)t},
\ee
that is (\ref{claim2}).\par
Let us now get the same type of estimates when $\lambda$ is not necessarily smaller than $\theta$. Let $\underline{u}_{\theta,0}$ be the function defined by:
$$\underline{u}_{\theta,0}(x)=\left\{\baa{ll}
\theta & \hbox{if }x\le-1,\\
-\theta\,x & \hbox{if }-1<x<0,\\
0 & \hbox{if }x\ge 0\eaa\right.$$
and denote by $\underline{u}_{\theta}$ the solution of the Cauchy problem (\ref{cauchy}) with initial condition $\underline{u}_{\theta,0}$. It follows then from (\ref{Omegatheta}) that
$$\forall\,T\ge\overline{t}_{\omega},\ \forall\,x\in\R,\quad u(T,x)\ge\underline{u}_{\theta,0}(x-y_{\omega}(T^+)),$$
whence
\be\label{uutheta}
\forall\,T\ge\overline{t}_{\omega},\ \forall\,t\ge 0,\ \forall\,x\in\R,\quad u(T+t,x)\ge\underline{u}_{\theta}(t,x-y_{\omega}(T^+))
\ee
from the maximum principle. But, as recalled in Section~\ref{intro}, there holds
$$\sup_{x\in\R}|\underline{u}_{\theta}(t,x)-\varphi_{c^*}(x-c^*t-m_{\theta}(t))|\to 0\ \hbox{ as }t\to+\infty,$$
where $m_{\theta}(t)=o(t)$ as $t\to+\infty$. In particular, given $\lambda\ge\theta$, there exists a time $T_{\lambda}$ (which also depends on $\theta$ and thus on $\epsilon$, but which does not depends on $T$) such that $\underline{u}_{\theta}(T_{\lambda},x)>\lambda$ for all $x<0$. Therefore,
$$\forall\,T\ge\overline{t}_{\omega},\ \forall\,x<y_{\omega}(T^+),\quad u(T+T_{\lambda},x)>\lambda$$
from (\ref{uutheta}). As a consequence, for all $t\ge\max(\overline{t}_{\omega}+T_{\lambda},t_{\lambda})$ and for all $x\in E_{\lambda}(t)$, one has $x\ge y_{\omega}((t-T_{\lambda})^+)$, whence
\be\label{claim2ter}
u_0(x)\le u_0(y_{\omega}((t-T_{\lambda})^+)=\omega\,e^{-\rho(t-T_{\lambda})}\le\Gamma \,e^{-(f'(0)-\epsilon)t}
\ee
for $t\ge\underline{T}_{\lambda,\epsilon,\Gamma}$, where $\underline{T}_{\lambda,\epsilon,\Gamma}\in[\max(\overline{t}_{\omega}+T_{\lambda},t_{\lambda}),+\infty)$ is large enough. This provides the claim (\ref{claim2}).\par
Estimates (\ref{claim1bis}), (\ref{claim2bis}) and (\ref{claim2ter}) yield (\ref{largetime}) with $T_{\lambda,\epsilon,\gamma,\Gamma}=\max(\overline{T}_{\lambda,\epsilon,\gamma},\underline{T}_{\lambda,\epsilon,\Gamma})$. That completes the proof of Theorem~\ref{th1}.\hfill$\Box$\break

Let us now turn to the proof of the more precise estimates of the location of the level sets of $u$ under the additional assumptions on $u_0$ and $f$ made in Theorem~\ref{th1bis}.\hfill\break

\noindent{\bf{Proof of Theorem~\ref{th1bis}.}} Let all assumptions of Theorem~\ref{th1} be fulfilled. Now, assume additionally that $u_0$ is convex in a neighborhood of $+\infty$, that is
$$u_0''(x)\ge 0\hbox{ for large }x.$$
We shall derive more precise lower bounds for $\min E_{\lambda}(t)$ at large time. Indeed, it follows from the arguments of the proof of part~c) of Theorem~\ref{th1} and from (\ref{subsolution}) that one can choose $\rho=f'(0)$, $\xi_2\ge\xi_0$ such that
$$|u_0''(x)|\le\frac{f'(0)\delta}{2(1+\delta)}\times u_0(x)\ \hbox{ for all }x\in[\xi_2,+\infty)$$
and $B=\max(s_1^{-\delta},2M/(f'(0)\delta))$, in such a way that the function $\underline{u}$ defined by (\ref{defunderu}) with $\rho=f'(0)$ is still a subsolution of (\ref{cauchy}). The last part of the proof then implies that, for every $\lambda\in(0,1)$, there exist a time $\underline{T}_{\lambda}\ge t_{\lambda}$ and a positive real number $\Gamma_{\lambda}$ such that
\be\label{Elambda}
E_{\lambda}(t)\subset u_0^{-1}\Big\{(0,\Gamma_{\lambda}\,e^{-f'(0)t}]\Big\}
\ee
for all $t\ge\underline{T}_{\lambda}$. Furthermore, it follows from (\ref{claim2bis}) with $\rho=f'(0)$ that the real numbers~$\Gamma_{\lambda}$ can be chosen independently of $\lambda$ when~$\lambda$ is small (under the notations of the proof of Theorem~\ref{th1}, one can choose $\Gamma_{\lambda}=\omega$ for $\lambda<\theta$).\par
Similarly, the same conclusion holds if there exists $\beta>0$ such that
$$u_0''(x)=O(u_0(x)^{1+\beta})\ \hbox{ as }x\to+\infty,$$
that is there exist $\xi'_0\ge\xi_0$ and $M'\ge 0$ such that $|u_0''(x)|\le M'\,u_0(x)^{1+\beta}$ in $[\xi'_0,+\infty)$. Indeed, one can choose $\beta'=\min(\beta,\delta)$, $\rho=f'(0)$, $\xi_2\ge\xi'_0$ such that
$$|u_0''(x)|\le\frac{f'(0)\beta'}{2(1+\beta')}\times u_0(x)\ \hbox{ for all }x\in [\xi_2,+\infty)$$
and $B=\max(s_1^{-\beta'},2(M+f'(0)M')/(f'(0)\beta'))$. Thus, the function $\underline{u}$ defined by~(\ref{defunderu}) with $\rho=f'(0)$ and $\beta'$ instead of $\delta$ is still a subsolution of (\ref{cauchy}): to check this point, one can first observe that (\ref{lowerf}) still holds with $\beta'$ instead of $\delta$, from (\ref{s10}) and since $0<\beta'\le\delta$. Then, in~(\ref{subsolution}), one can bound the term $-u_0''(x)\,e^{\rho t}$ by
$$|u_0''(x)\,e^{\rho t}|\le  M'\,u_0(x)^{1+\beta'}e^{\rho(1+\beta')t},$$
whence $\underline{u}_t(t,x)-\underline{u}_{xx}(t,x)-f(\underline{u}(t,x))\le 0$ provided that $\underline{u}(t,x)>0$. Eventually, the estimates (\ref{Elambda}) hold and the real numbers~$\Gamma_{\lambda}$ can be chosen independently of $\lambda$ when~$\lambda$ is small.\par
Lastly, consider the case when $f(s)\le f'(0)s-\mu\,s^{1+\nu}$ on $[0,1]$ and
$$|u_0''(x)|\le M'\,u_0(x)^{1+\beta}\le M'\,u_0(x)^{1+\nu}\ \hbox{ in }[\xi'_0,+\infty)$$
for some $\mu>0$, $\beta\ge\nu>0$, $M'\ge 0$ and $\xi'_0\ge\xi_0$. We shall derive upper bounds for the quantities $\max E_{\lambda}(t)$ for large $t$ which are more precise than the corresponding ones in~(\ref{largetime}) or~(\ref{largetimebis}). Set $\rho=f'(0)$, choose $\xi_1\ge\xi'_0$ such that
$$M'\,u_0(\xi_1)^{\nu}\le\mu,$$
and define $\overline{u}$ in $[0,+\infty)\times[\xi_1,+\infty)$ as in (\ref{defoveru}), that is
$$\overline{u}(t,x)=\min\left(\frac{u_0(x)\,e^{\rho t}}{u_0(\xi_1)},1\right).$$
As in the proof of Theorem~\ref{th1}, in order to check that $\overline{u}$ is a supersolution of (\ref{cauchy}) in $[0,+\infty)\times[\xi_1,+\infty)$, it remains to prove that $\overline{u}_t(t,x)-\overline{u}_{xx}(t,x)-f(\overline{u}(t,x))\ge 0$ as soon as $\overline{u}(t,x)<1$. For such a $(t,x)\in[0,+\infty)\times[\xi_1,+\infty)$ with $\overline{u}(t,x)<1$, there holds
$$\baa{rcl}
\overline{u}_t(t,x)-\overline{u}_{xx}(t,x)-f(\overline{u}(t,x)) & \ge & -\displaystyle{\frac{u_0''(x)\,e^{\rho t}}{u_0(\xi_1)}}+\displaystyle{\frac{\mu\,u_0(x)^{1+\nu}\,e^{\rho(1+\nu)t}}{u_0(\xi_1)^{1+\nu}}}\vspace{3pt}\\
& \ge & \left(\displaystyle{\frac{\mu}{u_0(\xi_1)^{\nu}}}-M'\right)\times\displaystyle{\frac{u_0(x)^{1+\nu}\,e^{\rho(1+\nu)t}}{u_0(\xi_1)}}\vspace{3pt}\\
& \ge & 0\eaa$$
due to the choice of $\xi_1$. Therefore, formulas (\ref{ineq1}), (\ref{ineq2}) and (\ref{ineq3}) hold with $\rho=f'(0)$. In particular, (\ref{ineq3}) says that, for any $\lambda\in(0,1)$, $t\ge t_{\lambda}$ and $y\in E_{\lambda}(t)$,
$$u_0(y)\ge\min(\eta,\lambda\,u_0(\xi_1)\,e^{-f'(0)t}),$$
where $\eta=\inf_{(-\infty,\xi_1)}u_0>0$. Thus, for every $\lambda\in(0,1)$, there exists a time $\overline{T}_{\lambda}\ge t_{\lambda}$ such that
$$\forall\, t\ge\overline{T}_{\lambda},\ \forall\, y\in E_{\lambda}(t),\quad u_0(y)\ge c\,\lambda\,e^{-f'(0)t},$$
where $c=u_0(\xi_1)>0$ does not depend on $\lambda$. This means (\ref{upperprecise}) and the proof of Theorem~\ref{th1bis} is complete.\hfill$\Box$\break

\noindent{\bf{Proof of Corollary~\ref{cor1}.}} Let $\xi:[0,+\infty)\to\R$ be any locally bounded function. From elementary arguments, it is straightforward to check that there exists a $C^1$ function $g:[0,+\infty)\to\R$ such that $g$ is of class $C^2$ on $[0,+\infty)\backslash E$, where $E$ is at most countable, and such that
$$\left\{\baa{l}
g(t)\ge\xi(2t)\ \hbox{ for all }t\ge 0,\vspace{3pt}\\
g'(t)>0\ \hbox{ for all }t\ge 0,\vspace{3pt}\\
g'(t)\to+\infty\hbox{ as }t\to+\infty,\vspace{3pt}\\
g''(t)=O(g'(t)^2)\hbox{ as }t\to+\infty\hbox{ in }[0,+\infty)\backslash E.\eaa\right.$$
Then, there exists a $C^1$ function $\widetilde{u}_0:\R\to(0,1)$ such that $\widetilde{u}'_0$ is bounded and negative on~$\R$ and
$$\widetilde{u}_0(x)=e^{-f'(0)\,g^{-1}(x)}\ \hbox{ for all }x\ge g(0)+1,$$
where $g^{-1}:[g(0),+\infty)\to[0,+\infty)$ denotes the reciprocal of the function $g$. In particular, the function $\widetilde{u}_0$ satisfies (\ref{frontlike}) and (\ref{slow}) since
$$\frac{\widetilde{u}'_0(x)}{\widetilde{u}_0(x)}=-\frac{f'(0)}{g'(g^{-1}(x))}\to 0\ \hbox{ as }x\to+\infty.$$
Furthermore,
$$\frac{\widetilde{u}''_0(x)}{\widetilde{u}_0(x)}=\frac{f'(0)\,g''(g^{-1}(x))}{(g'(g^{-1}(x)))^3}+\frac{f'(0)^2}{(g'(g^{-1}(x)))^2}\to 0\ \hbox{ as }x\to+\infty\hbox{ in }g([0,+\infty)\backslash E).$$
Let $\varphi$ be a nonnegative $C^{\infty}(\R)$ function whose support is included in $[-1,1]$ and whose integral over $\R$ is equal to $1$. Define $u_0=\varphi*\widetilde{u}_0(\cdot-1)$, that is
$$u_0(x)=\int_{\R}\varphi(y)\,\widetilde{u}_0(x-1-y)\,dy$$
for all $x\in\R$. The function $u_0$ is of class $C^{\infty}$, it ranges in $(0,1)$, it is decreasing with bounded derivative $u'_0$, and it satisfies (\ref{frontlike}), (\ref{slow}) and $u_0''(x)=o(u_0(x))$ as $x\to+\infty$. Furthermore, $u_0(x)\ge\widetilde{u}_0(x)$ for all $x\in\R$.\par
Let now $u$ be the solution of Cauchy with such an initial condition $u_0$, and let us prove that $u$ satisfies the conclusion of Corollary~\ref{cor1}. Fix any real number $\lambda\in(0,1)$. For all $t\ge t_{\lambda}$, the level set $E_{\lambda}(t)$ is a singleton $\{x_{\lambda}(t)\}$ (remember that $u_0$ is decreasing and $u(t,\cdot)$ is also decreasing for every $t\ge 0$). It follows from Theorem~\ref{th1} (and estimates (\ref{largetimebis}) applied with $\gamma=\Gamma=1$) that, for any $\epsilon\in(0,f'(0)/2)$, there exists a time $T_{\lambda,\epsilon}\ge t_{\lambda}$ such that
$$u_0^{-1}(e^{-(f'(0)-\epsilon)t})\le x_{\lambda}(t)\le u_0^{-1}(e^{-(f'(0)+\epsilon)t})$$
for all $t\ge T_{\lambda,\epsilon}$. Observe that $u_0^{-1}(e^{-(f'(0)-\epsilon)t})\ge\widetilde{u}_0^{-1}(e^{-(f'(0)-\epsilon)t})$. Therefore, one concludes that, for large $t$,
$$\min E_{\lambda}(t)=x_{\lambda}(t)\ge g\left(\frac{f'(0)-\epsilon}{f'(0)}\,t\right)\ge g\left(\frac{t}{2}\right)\ge\xi(t)$$
and the proof of Corollary~\ref{cor1} is complete.\hfill$\Box$


\SE{Uniform flatness at large time}\label{sec3}

This section is devoted to the proof of Theorem~\ref{th2}. Under the assumptions of this theorem, we will actually prove (\ref{flatbis}), which implies the flatness property (\ref{flat}), since $u$ stays bounded. The idea consists in using the special structure of the equation satisfied by the function~$u_x/u$ and in deriving some integral estimates which force this function to converge uniformly to~$0$ as $t\to+\infty$.\hfill\break

\noindent{\bf{Proof of Theorem~\ref{th2}.}} Set
$$v(t,x)=\frac{u_x(t,x)}{u(t,x)}\ \hbox{ for all }t\ge 0\hbox{ and }x\in\R.$$
We recall that $u_0>0$ in $\R$ and $u>0$ in $(0,+\infty)\times\R$ from the strong maximum principle. In Theorem~\ref{th2}, the function $v(0,\cdot)$ is assumed to be continuous, in $L^{\infty}(\R)\cap L^p(\R)\cap C^{2,\theta}(\R)$ (for some $1<p<+\infty$ and $0<\theta<1$). In particular, $v(0,x)\to 0$ as $x\to\pm\infty$. A direct computation shows that $v$ satisfies
\be\label{eqv}
v_t=v_{xx}+2\,v\,v_x+\left(f'(u)-\frac{f(u)}{u}\right)v,\quad t\ge 0,\ \ x\in\R.
\ee
Since $f'(s)\le f(s)/s$ for all $s\in(0,1)$, the maximum principle implies that $t\mapsto\|v(t,\cdot)\|_{L^{\infty}(\R)}$ is nonincreasing, and even that
\be\label{defMpm}\left\{\baa{l}
t\mapsto M^+(t):=\displaystyle{\mathop{\sup}_{x\in\R}}\ v(t,x)^+\hbox{ is nonincreasing},\vspace{3pt}\\
t\mapsto M^-(t):=\displaystyle{\mathop{\inf}_{x\in\R}}\ (-v(t,x)^-)\hbox{ is nondecreasing},\eaa\right.
\ee
where $v(t,x)^+=\max(v(t,x),0)$ and $v^-(t,x)=\max(-v(t,x),0)$.\par
Let us now check that $v(t,x)\to 0$ as $x\to\pm\infty$ for all $t\ge 0$. Choose an arbitrary $\epsilon>0$. Define $\alpha=2\,\|v(0,\cdot)\|_{L^{\infty}(\R)}$. Let $C>0$ be such that $\epsilon+C\,e^{-\alpha x}\ge v(0,x)$ for all $x\in\R$, and set
$$\overline{v}(t,x)=\epsilon+C\,e^{-\alpha(x-2\alpha t)}$$
for all $(t,x)\in[0,+\infty)\times\R$. Denote
$$b(t,x)=2\,v(t,x)\ \hbox{ and }\ c(t,x)=f'(u(t,x))-\frac{f(u(t,x))}{u(t,x)}\le 0$$
for all $t\ge 0$ and $x\in\R$. There holds $\overline{v}(0,x)\ge v(0,x)$ for all $x\in\R$, and
$$\overline{v}_t-\overline{v}_{xx}-b(t,x)\,\overline{v}_x-c(t,x)\,\overline{v}=C\,(\alpha^2+b(t,x)\alpha)e^{-\alpha(x-2\alpha t)}-c(t,x)\overline{v}(t,x)\ge 0$$
in $[0,+\infty)\times\R$, since $c\le 0$, $\overline{v}\ge 0$ and $\|b(t,\cdot)\|_{L^{\infty}(\R)}\le\alpha$ for all $t\ge 0$. The maximum principle yields $v\le\overline{v}$ in $[0,+\infty)\times\R$, whence
$$\limsup_{x\to+\infty}\,v(t,x)\le\epsilon\ \hbox{ for all }t\ge 0.$$
Similarly, one can prove that $\liminf_{x\to+\infty}v(t,x)\ge-\epsilon$ and $\limsup_{x\to-\infty}|v(t,x)|\le\epsilon$ for all $t\ge 0$. Furthermore, it follows from the choice of comparison functions that the limits are locally uniform in $t\ge 0$. In other words, $v(t,x)\to0$ as $x\to\pm\infty$ locally uniformly with respect to the variable $t\ge 0$.\par
Choose now an odd integer $q$ such that $p-1\le q$. Notice that $v(0,\cdot)\in L^{q+1}(\R)$. Let us check that $v(t,\cdot)\in L^{q+1}(\R)$ and $\|v(t,\cdot)\|_{L^{q+1}(\R)}\le\|v(0,\cdot)\|_{L^{q+1}(\R)}$ for all $t\ge 0$. Fix a time $t>0$. For all $R>0$, multiply (\ref{eqv}) by $v^q(t,x)$ and integrate by parts over $[0,t]\times[-R,R]$. It follows that
$$\baa{rcl}
\displaystyle{\frac{1}{q+1}}\displaystyle{\int_{-R}^R}(v^{q+1}(t,x)-v^{q+1}(0,x))\,dx & = & \displaystyle{\int_0^t}(v_x(s,R)\,v^q(s,R)-v_x(s,-R)\,v^q(s,-R))\,ds\vspace{3pt}\\
& & -q\displaystyle{\int\!\!\!\int_{[0,t]\times[-R,R]}}v_x^2(s,x)\,v^{q-1}(s,x)\,ds\,dx\vspace{3pt}\\
& & +\displaystyle{\frac{2}{q+2}}\displaystyle{\int_0^t}(v^{q+2}(s,R)-v^{q+2}(s,-R))\,ds\vspace{3pt}\\
& & +\displaystyle{\int\!\!\!\int_{[0,t]\times[-R,R]}}c(s,x)\,v^{q+1}(s,x)\,ds\,dx.\eaa$$
Since $c\le 0$ and $q\pm 1$ are even, one gets that
$$\baa{rcl}
\displaystyle{\frac{1}{q+1}}\displaystyle{\int_{-R}^R}(v^{q+1}(t,x)-v^{q+1}(0,x))\,dx & \le & \displaystyle{\int_0^t}(v_x(s,R)\,v^q(s,R)-v_x(s,-R)\,v^q(s,-R))\,ds\vspace{3pt}\\
& & +\displaystyle{\frac{2}{q+2}}\displaystyle{\int_0^t}(v^{q+2}(s,R)-v^{q+2}(s,-R))\,ds.\eaa$$
Since $v_x$ is bounded in $[0,t]\times\R$ and $v(s,x)\to 0$ as $x\to\pm\infty$ uniformly in $s\in[0,t]$, the right-hand side of the previous inequality converges to $0$ as $R\to+\infty$. As a consequence, $v(t,\cdot)\in L^{q+1}(\R)$ and $\|v(t,\cdot)\|_{L^{q+1}(\R)}\le\|v(0,\cdot)\|_{L^{q+1}(\R)}$.\par
In order to conclude that $\|v(t,\cdot)\|_{L^{\infty}(\R)}\to 0$ as $t\to+\infty$, one just needs to prove that $M^{\pm}(t)\to0$ as $t\to+\infty$, where $M^{\pm}(t)$ are defined in (\ref{defMpm}). We just do it for~$M^+(t)$, the case of $M^-(t)$ being similar. Remember that $t\mapsto M^+(t)$ is nonincreasing and nonnegative. Assume ab absurdo that $M^+(t)\to m>0$ as $t\to+\infty$. Then, there exists a sequence $(t_n,x_n)_{n\in\N}$ in $[0,+\infty)\times\R$ such that $t_n\to+\infty$ and $v(t_n,x_n)\to m$ as $n\to+\infty$. Define
$$v_n(t,x)=v(t+t_n,x+x_n)\ \hbox{ and }c_n(t,x)=c(t+t_n,x+x_n)$$
for all $n\in\N$ and $(t,x)\in[-t_n,+\infty)\times\R$. Remember that the sequences $(v_n)$ and $(c_n)$ are globally bounded, and that $c_n\le 0$ in $[-t_n,+\infty)\times\R$ for all $n\in\N$. Up to extraction of a subsequence, one can assume that $c_n\rightharpoonup c_{\infty}\le 0$ in $L^{\infty}_{loc}(\R\times\R)$ weak-* as $n\to+\infty$. From standard parabolic estimates, one can assume, up to extraction of another subsequence, that $v_n\to v_{\infty}$ in all $W^{1,p}_{loc}(\R\times\R)$ with respect to $t$ and $W^{2,p}_{loc}(\R\times\R)$ with respect to $x$, for all $1\le p<+\infty$, where $v_{\infty}$ is a bounded solution of
$$v_{\infty,t}=v_{\infty,xx}+2\,v_{\infty}\,v_{\infty,x}+c_{\infty}\,v_{\infty}$$
and $\max_{\R\times\R}v_{\infty}=v_{\infty}(0,0)=m>0$. Since $c_{\infty}\le 0$, one concludes from the strong maximum principle that $v_{\infty}(t,x)=m$ for all $(t,x)\in(-\infty,0]\times\R$. But since $v_n\to v_{\infty}$ (at least) locally uniformly in $\R\times\R$ as $n\to+\infty$, it follows that $\|v(t_n,\cdot)\|_{L^{q+1}(\R)}=\|v_n(0,\cdot)\|_{L^{q+1}(\R)}\to+\infty$ as $n\to+\infty$, which leads to a contradiction.\par
As a conclusion, $\|v(t,\cdot)\|_{L^{\infty}(\R)}\to 0$ as $t\to+\infty$. Since $u$ is globally bounded, this implies that $\|u_x(t,\cdot)\|_{L^{\infty}(\R)}\to 0$ as $t\to+\infty$ and the proof of Theorem~\ref{th2} is complete.\hfill$\Box$


\end{document}